\documentclass[a4paper,reqno]{amsart}
\usepackage{amssymb, amsmath, amscd}

\def\XIx\langle#1\rangle{h(#1)}
\newtheorem{theorem}{Theorem}[section]

\newtheorem{lemma}[theorem]{Lemma}

\newtheorem{remark}[theorem]{Remark}

\def\qedbox{\hbox{$\rlap{$\sqcap$}\sqcup$}}

\makeatletter
 \@addtoreset{equation}{section}
\makeatother
\begin{document}
\title[K\"ahler and para-K\"ahler curvature Weyl manifolds]{K\"ahler and para-K\"ahler curvature Weyl manifolds}
\author{Peter  Gilkey and Stana Nik\v cevi\'c}
\address{PG: Mathematics Department, \; University of Oregon, \;\;
  Eugene \; OR 97403 \; USA\\
  E-mail: gilkey@uoregon.edu}
\address{SN: Mathematical Institute, Sanu,
Knez Mihailova 36, p.p. 367,
11001 Belgrade,
Serbia.\\ Email: stanan@mi.sanu.ac.rs}
\begin{abstract}{We show that the Weyl structure of an almost-Hermitian Weyl manifold of dimension $n\ge6$ is trivial
if the associated curvature operator satisfies the K\"ahler identity. Similarly if the curvature of
an almost para-Hermitian Weyl manifold of dimension $n\ge6$ satisfies the para-K\"ahler identity,
then the Weyl structure is trivial as well.}\end{abstract}
\maketitle

\section{\bf Introduction}\label{sect-1}

\subsection{Pseudo-Riemannian Weyl geometry} Let $N$ be a smooth manifold of dimension $n\ge3$. Let $\nabla$ be a torsion free
connection on the tangent bundle $TN$ of $N$ and let $g$ be a pseudo-Riemannian metric on $N$ of signature $(p,q)$. Motivated by the
seminal paper of Weyl \cite {W22}, the triple
$\mathcal{W}:=(N,g,\nabla)$ is said to be a {\it Weyl manifold} if there exists a smooth $1$-form $\phi_{\nabla,g}\in C^\infty(T^*N)$ so
that:
\begin{equation}\label{eqn-1.a}
\nabla g=-2\phi_{\nabla,g}\otimes g\,.
\end{equation}

Weyl \cite{W22} used these geometries in an attempt to unify gravity
with electromagnetism -- although this approach failed for physical reasons, the resulting geometries are still an active area of investigation today. We
refer, for example, to \cite{G09} which studies Weyl geometry in the context of contact manifolds, to \cite{N09} where Einstein-Weyl structures
are examined in Lorentzian signature, to \cite{P10} where projectively flat Weyl manifolds are investigated, and to \cite{V10}
where the associated mass of an asymptotically flat Weyl structure is defined.

Let $[g]$ be the associated conformal class; $g_1\in[g]$ if and only
if there exists a smooth function $f$ so $g_1=e^{2f}g$. 
Weyl geometry is linked with conformal geometry as Equation (\ref{eqn-1.a}) means that $[g]$ is preserved by
covariant differentiation. If $g_1\in[g]$ and if $\mathcal{W}=(N,g,\nabla)$ is a Weyl manifold, then
the triple $\mathcal{W}_1:=(N, g_1,\nabla)$ is again a Weyl manifold where the associated
$1$-form is given by taking $\phi_{\nabla,g_1}:=\phi_{\nabla,g} - df$. We say the Weyl structure is {\it trivial} if there exists
$g_1\in[g]$ so that $\nabla=\nabla^{g_1}$ is the Levi-Civita connection of the metric $g_1$; additional equivalent conditions are given
below in Theorem \ref{thm-1.1}.

Let $\mathcal{R}$ be the {\it curvature operator} and let $R$ be  the
associated {\it curvature tensor} of the connection $\nabla$ of a Weyl manifold $\mathcal{W}=(N,g,\nabla)$:
\begin{eqnarray*}
&&\mathcal{R}(x,y):=\nabla_x\nabla_y-\nabla_y\nabla_x-\nabla_{[x,y]},\\
&&R(x,y,z,w):=g(\mathcal{R}(x,y)z,w)\,.
\end{eqnarray*}
Since $\nabla$ is torsion free, we have the symmetries:
\begin{eqnarray}
&&R(x,y,z,w)=-R(y,x,z,w),\label{eqn-1.b}\\
&&0=R(x,y,z,w)+R(y,z,x,w)+R(z,x,y,w)\label{eqn-1.c}\,.
\end{eqnarray}

The Ricci tensor is defined by setting:
\begin{equation}\label{eqn-1.d}
\operatorname{Ric}(x,y):=\operatorname{Tr}\{z\rightarrow \mathcal{R}(z,x)y\}\,.
\end{equation}
There is an additional well known curvature symmetry which pertains in Weyl geometry (see, for example, the discussion in \cite{GNS11}):
\begin{equation}\label{eqn-1.e}
R(x,y,z,w)+R(x,y,w,z)=\textstyle\frac2n\{\operatorname{Ric}(y,x) - \operatorname{Ric}(x,y)\}g(z,w)\,.
\end{equation}

If the Weyl structure is trivial, then $\nabla=\nabla^{g_1}$ for some $g_1\in[g]$ and we have the additional curvature symmetry
for the curvature $R^{g_1}$ of the Levi-Civita connection:
\begin{equation}\label{eqn-1.f}
R^{g_1}(x,y,z,w)+R^{g_1}(x,y,w,z)=0\,.
\end{equation}

We say that the curvature of $\mathcal{W}$ is {\it Riemannian} if in addition to the symmetries of Equation (\ref{eqn-1.b}) and of
Equation (\ref{eqn-1.c}), the symmetry of Equation (\ref{eqn-1.f}) is satisfied -- note that these 3 symmetries are conformal
invariants and that Equation (\ref{eqn-1.f}) implies Equation (\ref{eqn-1.e}). We have the following curvature condition which
ensures that the Weyl structure is trivial
\cite{GNS11}; we give the proof for the sake of completeness in Section \ref{sect-2.1}.

\begin{theorem}\label{thm-1.1}
Let $\mathcal{W}=(N,g,\nabla)$ be a Weyl manifold with $H^1(N;\mathbb{R})=0$. The following assertions are equivalent and if any is
satisfied, then the Weyl structure is trivial.
\begin{enumerate}
\item $d\phi_{\nabla,g}=0$.
\item $\nabla=\nabla^{g_1}$ for some $g_1\in[g]$.
\item $\nabla=\nabla^{g_1}$ for some pseudo-Riemannian metric $g_1$.
\item The curvature of $\nabla$ is Riemannian.
\end{enumerate}
\end{theorem}

\subsection{Almost para/pseudo-Hermitian Weyl geometry} Let $n=2\bar n\ge4$. We say that $(N,g,\nabla,J_-)$
is an {\it almost pseudo-Hermitian Weyl manifold} if $(N,g,\nabla)$ is a Weyl manifold, if $J_-$ is an {\it almost complex structure} on $TN$ (i.e.
$J_-$ is an endomorphism of $TN$ with
$J_-^2=-\operatorname{id}$), and if $J_-^*g=g$; necessarily $g$ has signature $(2\bar p,2\bar q)$ in this instance. Similarly, we say that
$(N,g,\nabla,J_+)$ is an {\it almost para-Hermitian Weyl manifold} if
$(N,g,\nabla)$ is a Weyl manifold, if
$J_+$ is a {\it para-complex structure} on $N$ (i.e. an endomorphism of
$TN$ with $J_+^2=\operatorname{id}$ and $\operatorname{Tr}(J_+)=0$), and if $J_+^*g=-g$; necessarily $g$ has neutral signature $(\bar n,\bar n)$.

The $\pm$ formalism permits us to discuss para-complex ($+$) and complex ($-$) geometry in parallel. For example, {\it (para)-Nijenhuis tensor} of an
almost (para)-complex manifold $(M,J_\pm)$ is given by
\begin{equation}\label{eqn-1.g}
N_\pm(x,y):=[x,y]\mp J_\pm [J_\pm x,y]\mp J_\pm [x,J_\pm y]\pm[J_\pm x,J_\pm y]\,.
\end{equation}
It vanishes if and only if $J_\pm$ is an integrable almost (para)-complex structure, i.e. given any point
$P\in N$, there exist local coordinates $(x^1,\dots,x^n)$ centered at $P$ so
\begin{equation}\label{eqn-1.h}
J_\pm\partial_{x_{2i-1}}=\partial_{x_{2i}}\quad\text{and}\quad
  J_\pm\partial_{x_{2i}}=\pm\partial_{x_{2i-1}}\quad\text{for}\quad 1\le i\le\bar n\,.
\end{equation}

\subsection{(Para)-K\"ahler Weyl geometry} Let
$(N,g,\nabla,J_\pm)$ be an almost para/pseudo-Hermitian Weyl manifold. If
$\nabla(J_\pm)=0$, then one says that this is a {\it (para)-K\"ahler Weyl manifold}. Note that necessarily
$J_\pm$ is integrable in this setting.
The study of such manifolds is very much an active research endeavor. See, for example, \cite{KK10} where the Siu-Beauville theorem
is extended to a certain class of compact K\"ahler-Weyl manifolds.

Pedersen, Poon, and Swann \cite{PPS93} used work of Vaisman \cite{V82,V83} to establish the following result in the Hermitian (i.e. positive definite)
setting; the extension to the higher signature setting and to the para-K\"ahler setting is immediate. We shall present their proof in Section
\ref{sect-2.3} for the sake of completeness.

\begin{theorem}\label{thm-1.2}
If $(N,g,\nabla,J_\pm)$ is a (para)-K\"ahler Weyl manifold with dimension $n\ge6$ and with $H^1(N;\mathbb{R})=0$, then the underlying Weyl structure
on $N$ is trivial.
\end{theorem}

We remark that Theorem \ref{thm-1.2} fails if $n=4$; see, for example, \cite{CP00,PS91}.

\subsection{Curvature (para)-K\"ahler Weyl manifolds} 
Suppose $(N,g,\nabla,J_\pm)$ is an almost para/pseudo-Hermitian Weyl manifold. If $\nabla(J_\pm)=0$, then one
has an additional curvature symmetry called the {\it K\"ahler identity}:
\begin{equation}\label{eqn-1.i}
\begin{array}{l}
\mathcal{R}(x,y)J_\pm=J_\pm\mathcal{R}(x,y)\quad\forall\quad x,y,\quad\text{or equivalently}\\
R(x,y,z,w)=\mp R(x,y,J_\pm z,J_\pm w)\quad\forall\quad x,y,z,w\,.
\end{array}\end{equation}
We say that $(N,g,\nabla,J_\pm)$ is a {\it (para)-K\"ahler
curvature Weyl manifold} if Equation (\ref{eqn-1.i}) is satisfied.
We will show in Section
\ref{sect-2.2} that there exist (para)-K\"ahler curvature Weyl manifolds where $J_\pm$ is not integrable; thus, in
particular, these are not (para)-K\"ahler Weyl manifolds.

The main result of this paper is the extension of Theorem
\ref{thm-1.2} to this context. The following result gives a curvature condition in these settings which
ensures that the Weyl structure is trivial; again it fails if $n=4$:
\begin{theorem}\label{thm-1.3}
If $(N,g,\nabla,J_\pm)$ is a curvature (para)-Ka\"hler Weyl manifold with dimension $n\ge6$ and with $H^1(N;\mathbb{R})=0$, then the underlying Weyl
structure on $N$ is trivial.
\end{theorem}

\subsection{Geometric realization results} It is convenient to work in a purely algebraic context. Let $V$ be a finite dimensional
vector space which is equipped with a non-degenerate symmetric bilinear form $h$ that we use to raise and
lower indices; the pair $(V,h)$ is said to be an {\it inner product space}. We say that
$A\in\otimes^4V^*$ is a {\it affine curvature tensor} if
$A$ has the symmetries given in Equations (\ref{eqn-1.b}) and (\ref{eqn-1.c}); let $\mathfrak{A}$ be the
set of all such tensors. The corresponding {\it affine curvature operator} $\mathcal{A}$ is defined by raising an index; $A$ and $\mathcal{A}$ are related
by the identity:
$$A(x,y,z,w)=\XIx\langle\mathcal{A}(x,y)z,w\rangle\quad\forall\quad x,y,z,w\in V\,.$$
 Let
$\mathfrak{W}$ be the subspace of
$\mathfrak{A}$ of all elements which in addition satisfy the symmetry of Equation (\ref{eqn-1.e}) and let $\mathfrak{R}$
be the subspace of
$\mathfrak{A}$ of elements which in addition satisfy the symmetry of Equation (\ref{eqn-1.f}); an element $A\in\mathfrak{R}$ is said to be a {\it
Riemannian curvature tensor} and the associated endomorphism $\mathcal{A}$ to be a {\it Riemannian curvature operator}. We have proper inclusions:
$$\mathfrak{R}\subset\mathfrak{W}\subset\mathfrak{A}\,.$$
The relations of Equations (\ref{eqn-1.b}) and (\ref{eqn-1.c}) generate the universal symmetries satisfied by the curvature of a
torsion free connection, the relations of Equations (\ref{eqn-1.b}), (\ref{eqn-1.c}), and (\ref{eqn-1.e}) generate the universal
symmetries satisfied by the curvature in Weyl geometry, and the relations of Equations (\ref{eqn-1.b}), (\ref{eqn-1.c}), and
Equation (\ref{eqn-1.f}) generate the universal symmetries satisfied by the curvature in pseudo-Riemannian geometry. We refer to
\cite{BGN11} for the proof of the following result and for other similar results and further bibliographic references concerning the
geometric realization of curvature in various contexts:

\begin{theorem}
Let $(V,h)$ be an inner product space.
\begin{enumerate}
\item If $A\in\mathfrak{A}$, then there exists a manifold $N$, there exists a point $P$ of $N$, there exists a
torsion free connection $\nabla$ on $TN$, and there exists an isomorphism $\Phi:T_PN\rightarrow V$ so that
$\Phi^*A=R_P$.
\item If $A\in\mathfrak{W}$, then there exists a Weyl manifold $(N,g,\nabla)$,  there exists a point $P$ of $N$, and there exists
an isomorphism
$\Phi:T_PN\rightarrow V$ so that $\Phi^*h=g_P$ and so that
$\Phi^*A=R_P$.
\item If $A\in\mathfrak{R}$, then there exists a pseudo-Riemannian manifold $(N,g)$, there exists a point $P$ of $N$, and
there exists an isomorphism $\Phi:T_N\rightarrow V$ so that $\Phi^*h=g_P$ and so that $\Phi^*A=R^g_P$.
\end{enumerate}
\end{theorem}

\subsection{Para/pseudo-Hermitian curvature models}
Let $(V,h)$ be an inner product space. We say that the triple $(V,h,J_\pm)$ is a {\it para/pseudo-Hermitian vector space} if $J_\pm$ is a (para)-complex
structure on
$V$ with $J_\pm^*h=\mp h$. Theorem
\ref{thm-1.3} will follow from Theorem
\ref{thm-1.1} and from the following purely algebraic result:

\begin{theorem}\label{thm-1.5}
Let $n\ge6$. Let $(V,h,J_\pm)$ be a para/pseudo-Hermitian vector space and let $A\in\mathfrak{W}$. If $A$ satisfies the (para)-K\"ahler
identity of {\rm Equation (\ref{eqn-1.i})}, then
$$A\in\mathfrak{R}\,.$$
\end{theorem}

Theorem \ref{thm-1.5} fails if $n=4$; there are non-trivial elements of $\mathfrak{W}-\mathfrak{R}$ which satisfy
the K\"ahler identity when $n=4$. We shall investigate this and related questions further in a subsequent paper.

\subsection{Outline of the paper}  In Section \ref{sect-2},
we prove Theorem \ref{thm-1.1}, we prove Theorem
\ref{thm-1.2}, and we exhibit a curvature (para)-K\"ahler manifold $(N,g,J_\pm)$ with (necessarily if $n\ge6$) trivial Weyl structure where $J_\pm$ is not
integrable.  In Section \ref{sect-3}, we review the basic group representation theory that we
shall need; these results are well known and we refer to the discussion in \cite{BGN11} Chapter 2 for example. We define the orthogonal group
$\mathcal{O}$, the (para)-unitary groups $\mathcal{U}_\pm$, and $\mathbb{Z}_2$ extensions $\mathcal{U}_\pm^*$ that play an important role in our
discussion. Suppose that
$G\in\{\mathcal{O},\mathcal{U},\mathcal{U}_\pm^*\}$. Results concerning the theory of submodules of
$\otimes^kV$ for the group $G$ are outlined in Section \ref{sect-3.2} and an introduction to the theory of scalar invariants
for $\otimes^kV$ is given in Section \ref{sect-3.3}. The para unitary group $\mathcal{U}_+$ is exceptional and these
results not apply to that group. 

In Section \ref{sect-4}, we review results of Singer and Thorpe \cite{ST69} decomposing $\mathfrak{R}$, results of Higa \cite{H93,H94} decomposing
$\mathfrak{W}$ as orthogonal modules, and an extension of results of Tricerri and Vanhecke \cite{TV81} decomposing $\mathfrak{R}$
as a $\mathcal{U}_\pm^\star$ module. These results are then used to decompose $\mathfrak{W}$ as a $\mathcal{U}_\pm^\star$ module. Theorem \ref{thm-1.5}
is then established Section
\ref{sect-5}. We refer to \cite{GI94, GNS11, H32, PS91, PT93} for further details concerning Weyl geometry.

\section{Geometric considerations}\label{sect-2}

\subsection{The proof of Theorem \ref{thm-1.1}}\label{sect-2.1}
Suppose that $d\phi_{\nabla,g}=0$. Since $H^1(N;\mathbb{R})=0$, we can express $\phi_{\nabla,g}=df$ for some function $f$. Let $g_1:=e^{2f}g\in[g]$.
Then $\phi_{\nabla,g_1}=0$ so $\nabla=\nabla^{g_1}$.
Thus Assertion (1) implies Assertion (2); by definition the Weyl structure is trivial if and only if Assertion (2) holds. Clearly Assertion (2)
implies Assertion (3). Since the curvature tensor of the Levi-Civita connection is Riemannian, Assertion (3) implies Assertion (4).
Suppose that Assertion (4) holds. We have $d\phi_{\nabla,g}=-\frac1n\Lambda\operatorname{Ric}$  where
$\Lambda\operatorname{Ric}$ is the alternating part of the Ricci tensor. Since the curvature tensor is Riemannian, the Ricci tensor is symmetric and
consequently
$\Lambda\operatorname{Ric}=0$. Thus Assertion (4) implies Assertion (1).\hfill\qed

\subsection{A curvature (para)-K\"ahler Weyl manifold which is not integrable}\label{sect-2.2}
Although relatively elementary, the following example is
instructive. Consider the usual coordinates $(x^1,\dots,x^n)$ on $N:=\mathbb{R}^n$. Let $J_-$ be
the standard complex structure given in Equation (\ref{eqn-1.h}). We work first in the positive definite setting. Let
\begin{equation}\label{eqn-2.a}
g:=dx^1\otimes dx^1+\dots+dx^{n}\otimes dx^n\,.
\end{equation}
Let $\Theta:\mathbb{R}^n\rightarrow \mathcal{O}$ satisfy
$\Theta(0)=\operatorname{id}$.
We consider a twisted almost complex structure:
$$J_-^\Theta:=\Theta^{-1}J_-\Theta\,.$$
Suppose that $\Theta=\Theta(x_1)$ is given by:
$$
\Theta\partial_{x_i}:=\left\{\begin{array}{lll}
\cos\theta(x_1)\partial_{x_1}+\sin\theta(x_1)\partial_{x_3}&\text{if}&i=1\\
\cos\theta(x_1)\partial_{x_3}-\sin\theta(x_1)\partial_{x_1}&\text{if}&i=3\\
\partial_{x_i}&\text{if}&i\ne1,3\end{array}\right\}.
$$
We compute the Nijenhuis tensor $N_-$ of Equation (\ref{eqn-1.g}) for this example. We have that $\{N_-(\partial_{x_1},\partial_{x_3})\}(0)$
consists of 4 parts:
\begin{enumerate}
\item $[\partial_{x_1},\partial_{x_3}](0)=0$.
\smallbreak\item $J_-^\Theta[J_-^\Theta\partial_{x_1},\partial_{x_3}](0)=-J_-(\partial_{x_3}J_-^\Theta)\partial_{x_1}=0$.
\smallbreak\item $J_-^\Theta[\partial_{x_1},J_-^\Theta\partial_{x_3}](0)
=\left.\left\{J_-(\partial_{x_1}J_-^\Theta)\partial_{x_3}\right\}\right|_{x=0}
   =\left.\left\{J_-\partial_{x_1}(\Theta^{-1}J_-\Theta)\partial_{x_3}\right\}\right|_{x=0}$
\smallbreak\quad $=\left.\left\{(-J_-\partial_{x_1}(\Theta) J_-+J_-J_-\partial_{x_1}(\Theta))\partial_{x_3}\right\}\right|_{x=0}$
\smallbreak\qquad $=\left.\left\{-J_-(\partial_{x_1}\Theta)|_{x=0}\partial_{x_4}-\partial_{x_1}(\Theta)|_{x=0}\partial_{x_3}\right\}\right.
=\partial_{x_1}|_{x=0}\ne0$.
\smallbreak\item
$-[J_-^\Theta\partial_{x_1},J_-^\Theta\partial_{x_3}](0)=-\left.\left\{(J_-\partial_{x_1})(J_-^\theta))\partial_{x_3}
-(J_-\partial_{x_3})(J_-^\Theta)\partial_{x_1}\right\}\right|_{x=0}$
\smallbreak\quad $=\left.\left\{(\partial_{x_2}(J_-^\theta)\partial_{x_3}-\partial_{x_4}(J_-^\theta)\partial_{x_1}\right\}\right|_{x=0}=0$.
\end{enumerate}
Thus the Nijenhuis tensor is non-trivial and $J_-^\Theta$ is not integrable. Since the curvature vanishes identically, $(N,g,J_-^\Theta)$ is
necessarily curvature K\"ahler. It is not, however, K\"ahler since $J_-^\theta$ is not integrable. By considering product manifolds, one can create
examples which are not flat. Furthermore, by replacing
$\cos$ and $\sin$ by
$\cosh$ and
$\sinh$ and modifying the signs appropriately, one can also construct examples in higher signature.

The construction of a curvature para-K\"ahler manifold which is not para-K\"ahler is similar. One replaces the complex structure
$J_-$ by the para-complex structure $J_+$ in Equation (\ref{eqn-1.h}), one replaces the
metric $g$ of Equation (\ref{eqn-2.a}) by the metric
$$g:=dx^1\otimes dx^1-dx^2\otimes dx^2+dx^3\otimes dx^3-dx^4\otimes dx^4\dots,$$
and one replaces the $N_-$ by $N_+$.
The remainder of the construction is unchanged and is therefore omitted.

\subsection{The proof of Theorem \ref{thm-1.2}}\label{sect-2.3} 
Let $(N,g,\nabla,J_\pm)$ be a (para)-K\"ahler Weyl manifold. Since $\nabla(J_\pm)=0$, $J_\pm$ is integrable. Let
$$\Omega_\pm(x,y):=g(x,J_\pm y)$$
be the associated {\it K\"ahler form}. We compute:
\begin{eqnarray*}
&&(\nabla_z\Omega_\pm)(x,y)=zg(x,J_\pm y)-g(\nabla_zx,J_\pm y)-g(x,J_\pm\nabla_zy)\\
&&\qquad=zg(x,J_\pm y)-g(\nabla_zx,J_\pm y)-g(x,\nabla_zJ_\pm y)\\
&&\qquad=(\nabla_zg)(x,J_\pm y)=-2\phi_{\nabla,g}(z)\Omega_\pm(x,y)\,.
\end{eqnarray*}
Let $\{e_i\}$ be a local frame for $TN$ and let $\{e^i\}$ be the dual frame for the cotangent bundle $T^*N$. We adopt the {\it
Einstein} convention and sum over repeated indices. Since
$\nabla$ is torsion free,
$d\Omega_\pm=e^i\wedge\nabla_{e_i}\Omega_\pm$.
Consequently
\begin{eqnarray*}
&&d\Omega_\pm=-2\phi_{\nabla,g}(e_i)e^i\wedge\Omega_\pm=-2\phi_{\nabla,g}\wedge\Omega_\pm,\\
&&0=d^2\Omega_\pm=-2d\phi_{\nabla,g}\wedge\Omega_\pm\,.
\end{eqnarray*}
Multiplication by $\Omega_\pm^{\frac12n-2}$ is an isomorphism between $\Lambda^2$ and $\Lambda^{n-2}$; this fact is usually cited
only in the positive definite setting for $J_-$ but extends to the more general situation. Thus as $n\ge6$,
$d\phi_{\nabla,g}\wedge\Omega_\pm=0$ implies $d\phi_{\nabla,g}=0$.\hfill\qedbox

\medbreak This argument fails if $n=4$; we can only conclude from this that $d\phi_{\nabla,g}\perp\Omega_\pm$.

\section{Representation theory}\label{sect-3}
In this section, we present the basic results from group representation theory that we shall need; these results are well known
and we refer, for example, to \cite{BGN11} Chapter 2 for further details. The structure groups are defined in Section
\ref{sect-3.1}. The theory of submodules of $\otimes^kV$ is outlined in Section \ref{sect-3.2}. Results relating to the theory of
scalar invariants are presented in Section \ref{sect-3.3}.
\subsection{Structure groups}\label{sect-3.1}
Let $h$ be a non-degenerate symmetric bilinear form on a real vector space $V$ of dimension $n$.
Let
$\mathcal{O}=\mathcal{O}(V,h)$ be the associated orthogonal group:
$$\mathcal{O}:=\{T\in\operatorname{GL}(V):T^*h=h\}\,.$$
If $(V,h,J_\pm)$ is a para/pseudo-Hermitian vector space, there are two associated Lie groups of interest. We define
the (para)-unitary group $\mathcal{U}_\pm$ and associated $\mathbb{Z}_2$ extension $\mathcal{U}_\pm^*$ by setting:
\begin{eqnarray*}
&&\mathcal{U}_\pm:=\{T\in\mathcal{O}:TJ_\pm=J_\pm T\},\\
&&\mathcal{U}_\pm^\star:=\{T\in\mathcal{O}:TJ_\pm=J_\pm T\text{ or }TJ_\pm=-J_\pm T\}\,.
\end{eqnarray*}

\subsection{Submodules of $\otimes^kV$}\label{sect-3.2}
We extend $h$ to $\otimes^kV$ so that
\begin{equation}\label{eqn-3.a}
\XIx\langle(v_1\otimes...\otimes v_k),(w_1\otimes...\otimes w_k)\rangle:=\prod_{i=1}^k\XIx\langle
v_i,w_i\rangle\,.
\end{equation}
Equation (\ref{eqn-3.a}) defines a non-degenerate symmetric bilinear form on $\otimes^kV$. We use $h$ to
identify
$V$ with $V^*$ and $\otimes^kV$ with $\otimes^kV^*$.
If $T\in\otimes^kV^*$ and if $T$ is a linear map of $V$, the {\it pull-back} $T^*\Theta$ is characterized by the identity
$$T^*\Theta(v_1,...,v_k)=\Theta(Tv_1,...,Tv_k)\,.$$
Let $G$ be one of the groups defined in Section \ref{sect-3.1}.  Then $G$ acts naturally on
$\otimes^kV^*$ by pull-back and preserves the canonical inner product defined in Equation (\ref{eqn-3.a}). Let $\xi$ be a $G$-invariant subspace of
$\otimes^kV^*$; the natural action of
$G$ on
$\otimes^kV^*$ makes $\xi$ into a $G$-submodule of $\otimes^kV$. The following is well known -- see, for example, the discussion in
\cite{BGN11} Chapter 2:
\begin{lemma}\label{lem-3.1}
Let $G\in\{\mathcal{O},\mathcal{U}_-,\mathcal{U}_\pm^*\}$. Let $\xi$ be a non-trivial $G$-submodule of $\otimes^kV^*$.
\begin{enumerate}\item $\xi$ is not totally isotropic.
\item There is an orthogonal direct sum decomposition $\xi=\eta_1\oplus...\oplus\eta_k$ where the $\eta_i$
are irreducible $G$-modules.
\item If $\xi_1$ and $\xi_2$ are inequivalent irreducible submodules of $\xi$, then $\xi_1\perp\xi_2$.
\item The multiplicity with which an irreducible representation appears in $\xi$ is independent of the decomposition in {\rm(1)}.
\item If $\xi_1$ appears with multiplicity $1$ in $\xi$ and if $\eta$ is any $G$-submodule of $\xi$, then either
$\xi_1\subset\eta$ or else $\xi_1\perp\eta$.
\end{enumerate}
\end{lemma}

\begin{remark}\label{rmk-3.2}
\rm
Much of what we will say subsequently extends to $\mathcal{U}_-$ with minor modifications. As the analysis of $\mathcal{U}_-$ is not needed to establish
the results of this paper, we shall not persue this topic. We note, however, that
Lemma \ref{lem-3.1} fails for the group $\mathcal{U}_+$. Let $(V,h,J_+)$ be a para-Hermitian
vector space. Decompose $V=V_+\oplus V_-$ into the $\pm1$ eigenspaces of $J_+$. Then $V_\pm$ are totally isotropic subspaces
of $V$ which are invariant under $\mathcal{U}_+$. 
\end{remark}

\subsection{Scalar invariants}\label{sect-3.3}
Let $\xi$ be a $G$-module. We say that $\Xi:\xi\rightarrow\mathbb{R}$ is a {\it scalar invariant} if $\Xi(g\cdot v)=\Xi(v)$ for
every $v\in\xi$ and for every $g\in G$; let $\mathcal{I}^G(\xi)$ be the vector space of all such invariants. Let
$\xi\subset\otimes^kV^*$.
H. Weyl \cite{W46} (see pages 53 and 66) gives a spanning set if $G=\mathcal{O}$ is the orthogonal group; the corresponding result
for the unitary group
$\mathcal{U}_-$ for in the Hermitian (i.e. positive definite) setting follows from
\cite{Fu58,Iw58} and the extension to the groups $\mathcal{U}_\pm^\star$ in general is straightforward -- see \cite{BGN11} for example. 

We discuss this spanning set. All invariants arise by using either the metric or the K\"ahler form to contract indices; invariants of
$\mathcal{U}_\pm^\star$ arise when the K\"ahler form appears an even number of times. It is worth being a bit more formal about this. Let
$(V,h,J_\pm)$ be a para/pseudo-Hermitian vector space. Let
$\Omega_{\pm,ij}$ be the components of the (para)-K\"ahler form. If $\{e_i\}$ is any basis for $V$, let $h_{ij}:=\XIx\langle e_i,e_j\rangle$.
The inverse matrix $h^{ij}=\XIx\langle e^i,e^j\rangle$ gives the components of the dual innerproduct  on $V^*$.
If $\Theta\in\otimes^{2k}V^*$, expand $\Theta=\Theta_{i_1\dots i_{2k}}e^{i_1}\otimes\dots\otimes e^{i_{2k}}$. Let $\pi\in\operatorname{Perm}(2k)$ be a
permutation of
$\{1,...,2k\}$. 
Let $\kappa_0:=h$, let $\kappa_1:=\Omega_\pm$, and let $\vec a$ be a sequence of $0$'s and $1$'s. Define:
$$\psi_{\pi,\vec a}(\Theta):=\kappa_{a_1}^{i_{\pi(1)}i_{\pi(2)}}
\dots\kappa_{a_k}^{i_{\pi(2k-1)}i_{\pi(2k)}}
\Theta_{i_1\dots i_{2k}}\,.$$
Let $n(\vec a)$ be the number of times $a_i=1$.
One then has:
\begin{lemma}\label{lem-3.3}
If $(V,h,J_\pm)$ is a para/pseudo-Hermitian vector space and if $\xi$ is a $\mathcal{U}_\pm^\star$ submodule of $\otimes^{2k}V^*$, then
$\mathcal{I}^{\mathcal{U}_\pm^\star}(\xi)=\operatorname{Span}_{n(\alpha)\text{ even}}\{\psi_{\pi,\alpha}\}$.
\end{lemma}

\section{Curvature decompositions}\label{sect-4}
In this section, we review the fundamental curvature decompositions that will play an important role our discussion. Section \ref{sect-4.1}
treats the Singer-Thorpe \cite{ST69} decomposition of $\mathfrak{R}$ as an $\mathcal{O}$ module. Section \ref{sect-4.2}
presents the Higa decomposition \cite{H93,H94} of $\mathfrak{W}$ as an $\mathcal{O}$ module. Section \ref{sect-4.3} discusses a
decomposition of $\mathfrak{R}$ as a $\mathcal{U}_\pm^\star$ module which generalizes the original Triceri-Vanhecke \cite{TV81} decomposition of
$\mathfrak{R}$ as a $\mathcal{U}_-$ module in the positive definite setting. Section
\ref{sect-4.4} gives the decomposition of
$\mathfrak{W}$ as a
$\mathcal{U}_\pm^\star$ module.

\subsection{The Singer-Thorpe  $\mathcal{O}$ module decomposition of $\mathfrak{R}$}\label{sect-4.1}
Let $\mathbb{R}\cdot  h\subset\otimes^2V^*$ be the trivial $1$-dimensional
$\mathcal{O}$ module, let $S_0^2\subset\otimes^2V^*$ be the
$\mathcal{O}$ module of trace free symmetric
$2$-tensors, and let $\Lambda^2\subset\otimes^2V^*$ be the $\mathcal{O}$ module of alternating $2$-tensors. Let
$P:=\ker\{\operatorname{Ric}\}\cap\mathfrak{R}$ be the $\mathcal{O}$ module of {\it Weyl conformal curvature tensors}. It
follows from \cite{ST69} that:
\begin{theorem} Let $n\ge4$. 
\begin{enumerate}
\item We may decompose $\otimes^2V^*=\mathbb{R}\cdot h\oplus S_0^2\oplus\Lambda^2$ as the orthogonal direct sum of
3 irreducible and inequivalent $\mathcal{O}$ modules.
\item There is an $\mathcal{O}$ isomorphism
$\mathfrak{R}\approx\mathbb{R}\oplus S_0^2\oplus P$ decomposing $\mathfrak{R}$ as the orthogonal direct sum of 3 irreducible and inequivalent
$\mathcal{O}$ modules.
\end{enumerate}
\end{theorem}

\subsection{The Higa $\mathcal{O}$ module decomposition of $\mathfrak{W}$}\label{sect-4.2}
If $\psi\in\Lambda^2$, define:
\begin{eqnarray}\label{eqn-4.a}
&&\sigma(\psi)(x,y,z,w):=2\psi(x,y)\XIx\langle z,w\rangle+\psi(x,z)\XIx\langle y,w\rangle-\psi(y,z)\XIx\langle x,w\rangle\\
&&\phantom{\sigma(\psi)(x,y,z,w):}-\psi(x,w)\XIx\langle y,z\rangle+\psi(y,w)\XIx\langle x,z\rangle\,.\nonumber
\end{eqnarray}
The map $\sigma$ is an $\mathcal{O}$ module isomorphism from $\Lambda^2$ to $\mathfrak{P}:=\sigma(\Lambda^2)$. We have \cite{H93,H94}:
\begin{theorem} Let $n\ge4$. We may decompose $\mathfrak{W}=\mathfrak{R}\oplus\mathfrak{P}$ as the orthogonal direct sum of $\mathcal{O}$ modules. This
gives a $\mathcal{O}$ module isomorphism
$\mathfrak{W}\approx\mathbb{R}\oplus S_0^2\oplus P\oplus\Lambda^2$ as the orthogonal direct sum of 4 irreducible and inequivalent $\mathcal{O}$ modules.
\end{theorem}

\subsection{The Tricerri-Vanhecke $\mathcal{U}_\pm^\star$ module decomposition}\label{sect-4.3}
The results of this section are the natural extension of results of Triceri and Vanhecke \cite{TV81} to the setting at hand and are discussed in
\cite{BGN11} in more detail. Let $(V,h,J_\pm)$ be a para/pseudo-Hermitian vector space. Define:
$$\begin{array}{ll}
S_{0,\mp}^{2,\mathcal{U}_\pm}:=\{\theta\in S^2:J_\pm^*\theta=\mp\theta\text{ and }\theta\perp h\},&
S_\pm^{2,\mathcal{U}_\pm}:=\{\theta\in S^2:J_\pm^*\theta=\pm\theta\},\\
\Lambda_{0,\mp}^{2,\mathcal{U}_\pm}:=\{\theta\in\Lambda^2:J_\pm^*\theta=\mp\theta\text{ and }\theta\perp\Omega_\pm\},&
\Lambda_\pm^{2,\mathcal{U}_\pm}:=\{\theta\in\Lambda^2:J_\pm^*\theta=\pm\theta\}.\vphantom{\vrule height 12pt}
\end{array}$$
\begin{lemma}\label{lem-4.1}
Let $n\ge4$. We may decompose
$$
   S^2=\mathbb{R}\cdot h\oplus S_{0,\mp}^{2,\mathcal{U}_\pm}\oplus S_\pm^{2,\mathcal{U}_\pm^\star}\text{ and }
  \Lambda^2=\mathbb{R}\cdot\Omega_\pm\oplus\Lambda_{0,\mp}^{2,\mathcal{U}_\pm^\star}\oplus\Lambda_\pm^{2,\mathcal{U}_\pm^\star}
$$
as the orthogonal direct sum of 6 irreducible and inequivalent $\mathcal{U}_\pm^\star$ modules.
\end{lemma}
\begin{remark}\rm The decomposition given above is also a decomposition of $S^2$ and $\Lambda^2$ into irreducible $\mathcal{U}_-$
modules. However $\mathbb{R}\cdot h\approx\mathbb{R}\cdot\Omega_+$ and $S_{0,+}^{2,\mathcal{U}_-}\approx \Lambda_{0,+}^{2,\mathcal{U}_-}$ as
$\mathcal{U}_-$ modules. In the para-Hermitian setting, we note that
$S_+^{2,\mathcal{U}_+}$ and $\Lambda_+^{2,\mathcal{U}_+}$ are not irreducible $\mathcal{U}_+$ modules.
\end{remark}

Let $n\ge8$. One follows \cite{TV81} to define $\mathcal{U}_-$ modules $ W_{-,i}$; these are also
$\mathcal{U}_-^*$ modules and there are analogous modules $\mathcal{U}_+^\star$ modules $ W_{+,i}$ in the para-Hermitian
setting. Set
$$\mathfrak{K}_{\pm,\mathfrak{R}}:=\{A\in\mathfrak{R}:A(x,y,z,w)=\mp A(x,y,J_\pm z,J_\pm w)\}\,;$$
these are the Riemannian curvature tensors which also satisfy the (para)-K\"ahler identity of Equation (\ref{eqn-1.i}).

\begin{theorem}\label{thm-4.5}
Let $(V,h,J_\pm)$ be a para/pseudo Hermitian vector space of dimension $n\ge8$. We may decompose
$$
  \mathfrak{R}=W_{\pm,1}\oplus...\oplus W_{\pm,10}\quad\text{and}\quad
  \mathfrak{K}_{\pm,\mathfrak{R}}=W_{\pm,1}\oplus  W_{\pm,2}\oplus  W_{\pm,3}
$$
as the orthogonal direct sum of irreducible
$\mathcal{U}_\pm^\star$ modules. We have
\begin{enumerate}\item $W_{\pm,1}\approx W_{\pm,4}\approx\mathbb{R}$ and
$ W_{\pm,2}\approx W_{\pm,5}\approx S_{0,\mp}^{2,\mathcal{U}_\pm}$.
\item 
$ W_{\pm,8}\approx S_\pm^{2,\mathcal{U}_\pm}$, and
$ W_{\pm,9}\approx\Lambda_\pm^{2,\mathcal{U}_\pm}$. 
\end{enumerate} With exception of the isomorphisms in (1), these are inequivalent
$\mathcal{U}_\pm$ modules.
\end{theorem}

\begin{remark}
\rm\ \begin{enumerate}
\item The original discussion of \cite{TV81} dealt with the unitary group $\mathcal{U}_-$ in the positive definite setting; we
refer to \cite{BGN11} for a discussion of the indefinite Hermitian setting and in the para-Hermitian setting.
If $n=6$, we set $ W_{\pm,6}=\{0\}$; if $n=4$, we set 
$ W_{\pm,5}= W_{\pm,6}= W_{\pm,10}=\{0\}$ to
achieve the corresponding decomposition. This does not affect our subsequent analysis.
\item Let $\Psi$ be the isomorphism from $\Lambda_\pm^2$ to
$ W_{\pm,9}$ given in (2) above; it is discribed quite explicitly in \cite{TV81} (page 372) in the Hermitian
setting and extends to our context to become:
\begin{eqnarray}
&&\Psi(\psi)(x,y,z,w):=2\XIx\langle x,J_\pm y\rangle\psi(z,J_\pm w)+2\XIx\langle z,J_\pm w\rangle\psi(x,J_\pm y)\nonumber\\
&&\qquad\qquad\qquad\qquad+\XIx\langle x,J_\pm z\rangle\psi(y,J_\pm w)+\XIx\langle y,J_\pm w\rangle\psi(x,J_\pm z)\label{eqn-4.b}\\
&&\qquad\qquad\qquad\qquad-\XIx\langle x,J_\pm w\rangle\psi(y,J_\pm z)-\XIx\langle y,J_\pm z\rangle\psi(x,J_\pm w)\,.\nonumber
\end{eqnarray}\end{enumerate}\end{remark}

\subsection{The decomposition of $\mathfrak{W}$ as a $\mathcal{U}_\pm^\star$ module}\label{sect-4.4}
Let $\sigma$ be as in Equation (\ref{eqn-4.a}).
We apply Lemma \ref{lem-4.1} to decompose $\Lambda^2$ and define:
$$W_{\pm,11}:=\sigma(\mathbb{R}\cdot\Omega_\pm),\quad
  W_{\pm,12}:=\sigma(\Lambda_{0,\mp}^{2,\mathcal{U}_\pm}),\quad
  W_{\pm,13}:=\sigma(\Lambda_\pm^{2,\mathcal{U}_\pm}).
$$
We combine Lemma \ref{lem-4.1} and Theorem \ref{thm-4.5} to establish:

\begin{theorem}\label{thm-4.7}
Let $(V,h,J_\pm)$ be a para/pseudo Hermitian vector space of dimension $n\ge8$. We may decompose
$$\mathfrak{W}=W_{\pm,1}\oplus...\oplus W_{\pm,13}$$
as the orthogonal direct sum of irreducible
$\mathcal{U}_\pm^\star$ modules. With the exception of the isomorphisms noted in {\rm Theorem \ref{thm-4.5}}, these are inequivalent
$\mathcal{U}_\pm^\star$ modules.
\end{theorem}

\begin{remark}
\rm As before, we shall set $ W_{\pm,6}=\{0\}$ if $n=6$ and we shall set
$ W_{\pm,5}= W_{\pm,6}= W_{\pm,10}=\{0\}$ if $n=4$.
The modules $\{\mathbb{R},S_{0,\mp}^{2,\mathcal{U}_\pm},\Lambda_\pm^{2,\mathcal{U}_\pm}\}$ appear with multiplicity 2 in the decomposition of
$\mathfrak{W}$ as a $\mathcal{U}_\pm^*$ module; the remaining modules appear with multiplicity 1.
\end{remark}

\subsection{The modules $\Lambda_\pm^{2,\mathcal{U}_\pm^\star}$}
We shall
need the following technical result:

\begin{lemma}\label{lem-4.9}
If $\xi$ is a non-trivial proper $\mathcal{U}_\pm^*$ submodule of $\Lambda_\pm^{2,\mathcal{U}_\pm}\oplus\Lambda_\pm^{2,\mathcal{U}_\pm}$, then there
exists
$(a,b)\ne(0,0)$ so 
$$\xi=\xi(a,b):=\{(a\theta,b\theta)\}_{\theta\in\Lambda_\pm^{2,\mathcal{U}_\pm}}
\subset\Lambda_\pm^{2,\mathcal{U}_\pm}\oplus\Lambda_\pm^{2,\mathcal{U}_\pm}\,.$$
\end{lemma}

\begin{proof} We have
\begin{eqnarray*}
&&\Lambda_\pm^{2,\mathcal{U}_\pm}\otimes\Lambda_\pm^{2,\mathcal{U}_\pm}
=\{\theta\in\otimes^4V^*:\theta(x,y,z,w)=-\theta(y,x,z,w)=-\theta(x,y,w,z)\\
&&\quad\text{and }\theta(x,y,z,w)=\pm\theta(J_\pm x,J_\pm y,z,w)=\pm\theta(x,y,J_\pm z,J_\pm w)\}\,.
\end{eqnarray*}
It follows from these symmetries and from Lemma \ref{lem-3.3} that there is only one $\mathcal{U}_\pm^\star$ invariant of
$\Lambda_\pm^{2,\mathcal{U}_\pm}\otimes\Lambda_\pm^{2,\mathcal{U}_\pm}$ given by $h^{ik}h^{jl}\theta(e_i,e_j,e_k,e_l)$. Thus
\begin{equation}\label{eqn-4.c}
\dim\{\mathcal{I}^{\mathcal{U}_\pm^*}(\Lambda_\pm^{2,\mathcal{U}_\pm}\otimes\Lambda_\pm^{2,\mathcal{U}_\pm})\}\le1\,.
\end{equation}
Let $\operatorname{Hom}^{\mathcal{U}_\pm^*}(\Lambda_\pm^{2,\mathcal{U}_\pm})$ be the set of all linear maps
$T:\Lambda_\pm^{2,\mathcal{U}_\pm}\rightarrow\Lambda_\pm^{2,\mathcal{U}_\pm}$ with $Tg=gT$ for all $g\in\mathcal{U}_\pm^*$. Let
$\Xi_T(\theta_1\otimes\theta_2):=\XIx\langle\theta_1,T\theta_2\rangle$ be the linear invariant defined by
$T\in\operatorname{Hom}^{\mathcal{U}_\pm^*}(\Lambda_\pm^{2,\mathcal{U}_\pm})$;
for example, $\Xi_1=\Xi_{\operatorname{id}}$. Equation (\ref{eqn-4.c}) then shows
\begin{equation}\label{eqn-4.d}
\operatorname{Hom}^{\mathcal{U}_\pm^*}(\Lambda_\pm^{2,\mathcal{U}_\pm^\star})=\operatorname{Id}\cdot\mathbb{R}\,.
\end{equation}
Let $\xi$ be a proper $\mathcal{U}_\pm^*$ submodule of
$\Lambda_\pm^{2,\mathcal{U}_\pm}\oplus\Lambda_\pm^{2,\mathcal{U}_\pm}$. Let $\pi_1$ (resp. $\pi_2$) be
projection on the first (resp. on the second) factor. Since $\xi$ is non-trivial, we may assume without loss of generality that
$\pi_1\xi\ne\{0\}$; since $\xi$ is a proper submodule, $\xi$ is necessarily irreducible and hence $\pi_1$ is an isomorphism. If
$\pi_2=0$, then $\xi=\xi(1,0)$. Thus we may assume that $\pi_2\ne0$ and hence $\pi_2^{-1}\pi_1=T$ is a non-trivial $\mathcal{U}_\pm^*$
equivariant map of $\Lambda_\pm^{2,\mathcal{U}_\pm}$. Equation (\ref{eqn-4.d}) then shows $T=b\operatorname{id}$ and $\xi=\xi(1,b)$.
\end{proof}

\section{The proof of Theorem \ref{thm-1.5}}\label{sect-5}

Let $(V,h,J_\pm)$ be a para/pseudo-Hermitian vector space. Let
$$\mathfrak{K}_{\pm,\mathfrak{W}}:=\{A\in\mathfrak{W}:A(x,y,z,w)=\mp A(x,y,J_\pm z,J_\pm w)\ \forall\ x,y,z,w\}$$
be the space of all Weyl tensors satisfying the (para)-K\"ahler identity of Equation (\ref{eqn-1.i}). We use the decomposition of Theorem \ref{thm-4.7}
and set
$$
   \mathfrak{K}_{\pm,\mathfrak{W}}^1:=\left\{\oplus_{4\le i\le 13} W_{\pm,i}\right\}\cap\mathfrak{K}_{\pm,\mathfrak{W}}\,.
$$
We use Theorem \ref{thm-4.5} to see 
$\mathfrak{R}\cap\mathfrak{K}_{\pm,\mathfrak{W}}=W_{\pm,1}\oplus W_{\pm,2}\oplus W_{\pm,3}$. Consequently
$$
  \mathfrak{K}_{\pm,\mathfrak{W}}= W_{\pm,1}\oplus W_{\pm,2}
    \oplus W_{\pm,3}\oplus\mathfrak{K}_{\pm,\mathfrak{W}}^1\,.
$$

We prove Theorem \ref{thm-1.5} by showing $\mathfrak{K}_{\pm,\mathfrak{W}}^1=\{0\}$. Suppose that
$4\le i\le 13$ and
$i\ne 9,13$. Since $W_{\pm,i}$ appears with multiplicity $1$ in $\mathfrak{K}_{\pm,\mathfrak{W}}^1$, Lemma \ref{lem-3.1} shows that either
$W_{\pm,i}\subset\mathfrak{K}_{\pm,\mathfrak{W}}^1$ or
$W_{\pm,i}\perp\mathfrak{K}_{\pm,\mathfrak{W}}^1$.  By
Theorem \ref{thm-4.5},
$W_{\pm,i}\cap\mathfrak{K}_{\pm,\mathfrak{R}}=\{0\}$ for $4\le i\le 10$. Consequently
$$
   \mathfrak{K}_{\pm,\mathfrak{W}}^1=\left\{W_{\pm,9}\oplus W_{\pm,11}\oplus W_{\pm,12}\oplus W_{\pm,13}\right\}
   \cap\mathfrak{K}_{\pm,\mathfrak{W}}\,.
$$
\subsection{The module $W_{\pm,11}$}
We use Equation (\ref{eqn-4.a}) to see:
\begin{enumerate}
\item $\sigma(\Omega_\pm)(e_1,e_4,e_3,e_1)=-h(e_4,J_\pm e_3)h(e_1,e_1)=-h_{11}h_{44}$,
\smallbreak\item $\mp\sigma(\Omega_\pm)(e_1,e_4,J_\pm e_3,J_\pm e_1)=\pm h(e_1,J_\pm J_\pm e_1)h(e_4,J_\pm e_3)=h_{11}h_{44}$,
\smallbreak\item Thus $\sigma(\mathbb{R}\cdot\Omega_\pm)\not\subset\mathfrak{K}_{\pm,\mathfrak{W}}^1$ if $n\ge4$.
\end{enumerate}
\subsection{The module $W_{\pm,12}$}
Let $\psi_{0,\pm}:=e^1\otimes e^2-e^2\otimes e^1+\delta_\pm\{e^3\otimes e^4-e^4\otimes e^3\}$ where $\delta_\pm$ is chosen to ensure
$\psi_{0,\pm}\perp\Omega_\pm$. We have $J_\pm^*\psi_{0,\pm}=\mp\psi_{0,\pm}$ and thus $\psi_{0,\pm}\in\Lambda_{0,\pm}^{2,\mathcal{U}_\pm}$. We use
Equation (\ref{eqn-4.a}) to verify:
\begin{enumerate}
\item 
$\sigma(\psi_{0,\pm})(e_5,e_1,e_2,e_5)=-\psi_{0,\pm}(e_1,e_2)\XIx\langle e_5,e_5\rangle=-h_{55}$.
\smallbreak\item
$\mp\sigma(\psi_{0,\pm})(e_5,e_1,J_\pm e_2,J_\pm e_5)=\pm\psi_{0,\pm}(e_5,J_\pm e_5)\XIx\langle e_1,J_\pm e_2\rangle=0$.
\smallbreak\item
$W_{\pm,12}\not\subset\mathfrak{K}_{\pm,\mathfrak{W}}^1$ if $n\ge6$.
\end{enumerate}

\subsection{The module $W_{\pm,9}\oplus W_{\pm,13}$}
Let $\psi_\pm:=e^1\otimes e^3-e^3\otimes e^1\pm e^2\otimes e^4\mp e^4\otimes e^2$. 
Then $J_\pm^*\psi_\pm=\pm\psi_\pm$ so $\psi_\pm\in\Lambda_\pm^{2,\mathcal{U}_\pm}$. By Equation (\ref{eqn-4.a}) and Equation
(\ref{eqn-4.b}):
\begin{enumerate}
\item $\sigma(\psi_\pm)(e_5,e_1,e_3,e_5)=-\psi_\pm(e_1,e_3)\XIx\langle e_5,e_5\rangle=-h_{55}$.
\item $\sigma(\psi_\pm)(e_5,e_1,e_4,e_6)=0$.
\item $\Psi(\psi_\pm)(e_5,e_1,e_3,e_5)=0$.
\item $\Psi(\psi_\pm)(e_5,e_1,e_4,e_6)=-\psi_\pm(e_1,J_\pm e_4)\XIx\langle e_5,J_\pm e_6\rangle=-h_{55}$.
\item $\sigma(\psi_\pm)(e_5,e_6,e_1,e_4)=0$.
\item $\sigma(\psi_\pm)(e_5,e_6,J_\pm e_1,J_\pm e_4)=0$.
\item $\Psi(\psi_\pm)(e_5,e_6,e_1,e_4)=2\XIx\langle e_5,J_\pm e_6\rangle\psi_\pm(e_1,J_\pm e_4)=2h_{55}$.
\item $\Psi(\psi_\pm)(e_5,e_6,J_\pm e_1,J_\pm e_4)=2\XIx\langle e_5,J_\pm e_6\rangle\psi_\pm(J_\pm e_1,J_\pm J_\pm e_4)=\pm2h_{55}$.
\end{enumerate}
          
\medbreak For $(a,b)\ne(0,0)$, let $\xi(a,b):=\operatorname{Range}\{a\sigma+b\Psi\}\subset W_{\pm,9}\oplus W_{\pm,13}$. 
We suppose $\xi(a,b)\cap\mathfrak{K}_{\pm,\mathfrak{W}}^1\ne\{0\}$ and thus $\xi(a,b)\subset\mathfrak{K}_{\pm,\mathfrak{W}}^1$. Assertions (1)-(4)
then yield $a=\mp b$ while Assertions (5)-(8) yield $b=0$. We
apply Lemma \ref{lem-4.9} to see that every non-trivial proper submodule of
$W_{\pm,9}\oplus W_{\pm,13}$ is isomorphic to $\xi(a,b)$ for some $(a,b)\ne0$. Thus
$$\left\{W_{\pm,9}\oplus W_{\pm,13}\right\}\cap\mathfrak{K}_{\pm,\mathfrak{W}}^1=\{0\}\,.$$
and consequently $\mathfrak{K}_{\pm,\mathfrak{W}}^1=\{0\}$. This completes the proof of Theorem \ref{thm-1.5}.\hfill\qedbox

\section*{Acknowledgments} 
Research of P. Gilkey partially supported by DFG PI 158/4-6 (Germany) and by project MTM2009-07756 (Spain). Research of
S. Nik\v cevi\'c partially supported by a research grant of the TU Berlin, by project MTM2009-07756 (Spain), and by 144032
(Serbia).  It is a pleasant task to acknowledge useful conversations with E. Garcia-Rio of the Universidad de Santiago (Spain), with U. Simon of the
Technische Universit\"at (Berlin), and with A. Swann of the University of Southern Denmark.

\end{document}